%% file: mainArxiv.tex
\documentclass{amsart}

\usepackage[margin=1.4in]{geometry}

\input{ARXIVpackages}
\input{ARXIVdefinitions}

\newcommand{\fd}{.}

\title{Polyhedral Geometry in Oscar}
\author{Taylor Brysiewicz}
\address{Department of Mathematics, University of Western Ontario, London, Canada.}
\email{tbrysiew@uwo.ca}

\author{Michael Joswig}
\address{TU Berlin, Chair of Discrete Mathematics/Geometry, Berlin, Germany; Max-Planck-Institut für Mathematik in den Naturwissenschaften, Leipzig, Germany.}
\email{joswig@math.tu-berlin.de}

\usepackage{cite}
 
\begin{document}
\begin{abstract}
\OSCAR \cite{OSCAR} is an innovative new computer algebra system which combines and extends the power of its four cornerstone systems - GAP (group theory), Singular (algebra and algebraic geometry), Polymake (polyhedral geometry), and Antic (number theory). Assuming little familiarity with the subject, we give an introduction to computations in polyhedral geometry using \OSCAR, as a chapter of the upcoming \OSCAR book \cite{OSCAR-book}.
  In particular, we define polytopes, polyhedra, and polyhedral fans, and we give a brief overview about computing convex hulls and solving linear programs.
  Three detailed case studies are left for experts in polyhedral geometry.
  These are concerned with face numbers of random polytopes, constructions and properties of Gelfand--Tsetlin polytopes, and secondary polytopes.
\end{abstract}

\maketitle

\section{Introduction}
\label{ch:cs-polyhedral-geometry}
Polyhedral geometry is linear algebra over ordered fields.
The interest in the topic comes from several applications.
Polyhedra arise as feasible regions of linear programs, whence they occur naturally in linear and integer programming \cite{Schrijver:TOLIP}.
 Newton polytopes of multivariate polynomials arise, implicitly or explicitly, in algorithms for solving systems of polynomial equations \cite{Sturmfels:1996}; this has further consequences for toric \cite{CLS:2011} and tropical geometry \cite{MS:2015,ETC}.
There are further connections to commutative algebra, combinatorial topology and more; e.g., see \cite{Stanley:1996}.
This wide range of applications lead to polyhedral geometry becoming a field of its own \cite{Gruenbaum:2003,Ziegler:1995,DLRS:2010}.
The latter textbooks are recommended for further reading.
Our conventions follow \cite{JT:2013} and \cite{Ziegler:1995}.

From an algorithmic point of view, polyhedra are rather benign.
Since linear algebra suffices for most computations, we immediately get exact algorithms over any (ordered) field for which we have an exact implementation of the arithmetic.
Of course, the rational numbers form the prime examples, but real number fields \cite[p.94f]{Gruenbaum:2003} and fields of rational functions \cite{JoswigLohoLorenzSchroeter:2016} occur in practice, too.

\section{Polytopes and Polyhedra}
A \emph{polyhedron} in $\RR^d$ is the intersection of finitely many affine closed half-spaces. Equivalently, polyhedra can be described by finitely many linear inequalities, whence they are precisely the feasible regions of linear programs. A \emph{polytope} is a bounded polyhedron.

An (affine) hyperplane $H$ \emph{supports} a polyhedron $P$ if $H$ intersects $P$ nontrivially, but $H$ does not separate $P$, i.e., one of the two half-spaces defined by $H$ contains all of $P$.
The intersection of $P$ with a supporting hyperplane is a \emph{(proper) face} of $P$.
Each polyhedron has finitely many faces.
They are partially ordered by inclusion and by construction, are polyhedra themselves.
Each polyhedron has a \emph{dimension} $r$, defined as the dimension of its affine hull, and may be described as an $r$-polyhedron or $r$-polytope accordingly.
We allow for polyhedra which are not full dimensional; in that case we have $r<d$.

Furthermore, we also allow the empty set as a polytope.
In that case, there are no supporting hyperplanes nor proper faces.
We use the convention $\dim\emptyset=-1$.

\subsection{Construction and Basic Properties}

The convex hull of finitely many points is a polytope and conversely; see \cite[\S1.1]{Ziegler:1995} or \cite[\S3.1.3]{JT:2013}.
Let us see this in a first example. We create a polytope $P$ which is the convex hull of five points in the plane.
\begin{minted}{jlcon}
julia> points = [1 0; 1 1; 0 1; -1 0; -1 -1];

julia> P = convex_hull(points)
Polyhedron in ambient dimension 2
\end{minted}
Mathematically a lone pentagon is not so exciting.
Here it serves to demonstrate how to call the most basic polyhedral functions in \OSCAR, and we can see the full output.
The \emph{facets} of a polyhedron are its maximal proper faces. They yield a representation of the polyhedron in terms of a minimal number of linear inequalities.
If a polyhedron is full-dimensional, then its minimal inequality description is unique.
\begin{minted}{jlcon}
julia> facets(P)
5-element SubObjectIterator{AffineHalfspace{QQFieldElem}} over the Halfspaces of R^2 described by:
-x_1 <= 1
-x_1 + x_2 <= 1
x_1 - 2*x_2 <= 1
x_1 <= 1
x_2 <= 1
\end{minted}
The \emph{$f$-vector} counts the faces per dimension.
For the pentagon, we have five $0$-dimensional faces, the \emph{vertices} (which were our input), and the five $1$-dimensional facets that we saw just above.
\begin{minted}{jlcon}
julia> f_vector(P)
2-element Vector{ZZRingElem}:
 5
 5

julia> volume(P)
5//2
\end{minted}
As a compact subset of $\RR^d$ any polytope has a finite Euclidean volume.
The (Euclidean) volume in this planar example is the area.
The points with integer coordinates in a polyhedron are its \emph{lattice points}.
In our example the lattice points of $P$ are the five vertices along with the origin.
\begin{minted}{jlcon}
julia> lattice_points(P)
6-element SubObjectIterator{PointVector{ZZRingElem}}:
 [-1, -1]
 [-1, 0]
 [0, 0]
 [0, 1]
 [1, 0]
 [1, 1]
\end{minted}

Note that all coefficients and metric values we computed so far are exact rational numbers.
More generally, whenever a polytope is given by inequalities whose coefficients lie in some ordered field $\FF$, then the vertex coordinates also lie in $\FF$.
Via triangulation the volume is a sum of determinant expressions in the vertex coordinates, hence this lies in $\FF$ too.
\OSCAR allows coordinates in an arbitrary real number field; the corresponding coordinate type is \mintinline{jl}{EmbeddedAbsSimpleNumFieldElem}.

\subsection{Algorithmic Ingredients}
Going back to our initial pentagon example, we saw that a polytope which is defined as the convex hull of finitely points, admits an inequality description in terms of its facets.
The algorithmic problem of passing from the points to the inequalities is known as the \emph{convex hull problem}.
Several different methods are known and implemented, but the precise complexity status of the general convex hull problem is open \cite[Open Problem 26.3.4]{HDCG3:convex-hulls}.
In polyhedral geometry practice, computing convex hulls is usually a crucial step and often the bottleneck.
The algorithms and their implementations compared in the survey \cite{polymake:2017} are available in \polymake and thus also in \OSCAR.

Computing face lattices and $f$-vectors of polytopes starts with the incidences between the vertices and the facets, which are usually determined after (or simultaneously with) a convex hull computation.
State of the art is an algorithm of Kaibel and Pfetsch \cite{KaibelPfetsch:2002}, which allows for numerous variations \cite{HampeJoswigSchroeter:MEGA2017}.
Kliem and Stump describe a fast algorithm for computing the $f$-vector from the vertex-facet incidences without constructing the covering relations of the entire face lattice \cite{KliemStump:2022}.
For volume computation see Büeler, Enge and Fukuda \cite{BuelerEngeFukuda:2000}, for further details on several algorithms for polytopes see Kaibel and Pfetsch \cite{KaibelPfetsch:2003}.

\subsection{Linear Programs}
Polytopes and polyhedra occur most naturally in linear optimization.
The task of a \emph{linear program} is to optimize an affine linear function $x \mapsto c^Tx + k$, called the \emph{objective function}, over the points $x$ within some polyhedron $P \subseteq \RR^d$ called the \emph{feasible set}. 
For simplicity, we show examples just in the plane.
\begin{minted}{jlcon}
julia> P = polyhedron([1 5; 2 -1; -1 0; 0 -1], [20,10,0,0])
Polyhedron in ambient dimension 2

julia> LP = linear_program(P, [1,1], convention=:max)
Linear program
   max{c*x + k | x in P}
where P is a Polyhedron{QQFieldElem} and
   c=Polymake.LibPolymake.Rational[1 1]
   k=0

julia> optimal_value(LP)
100/11

julia> optimal_vertex(LP)
2-element PointVector{QQFieldElem}:
 70//11
 30//11
\end{minted}
The data of a linear program consists of a polyhedron, an objective function, and a choice of \mintinline{jl}{max} or \mintinline{jl}{min}.
Since the value of $k$ in the objective function does not influence the optimal solutions, it is zero by default. In our example, the maximum is attained at a unique point, namely the vertex $(70/11, 30/11)$.

If the objective function is bounded, with optimum, say, $\mu$, then the set $\smallSetOf{x\in P}{c^Tx +k = \mu}$ of optimal solutions forms a face of the feasible region $P$.
Conversely, each face of $P$ arises in this way.

\subsection{The Wonderful World of Convex Polytopes}
Convex polytopes occur throughout mathematics.

\subsubsection{The Platonic Solids and Their Friends:}
A $d$-polytope $P$ is \emph{regular} if its group of symmetries acts transitively on the flags of faces of $P$.
Specializing to $d=3$ this leads us to the \emph{Platonic solids}, the most classical polytopes of all.
Up to scaling (and rigid motions) there are exactly five Platonic solids.
The simplest one is the regular tetrahedron.
\begin{minted}{jlcon}
julia> T = tetrahedron()
Polytope in ambient dimension 3

julia> vertices(T)
4-element SubObjectIterator{PointVector{QQFieldElem}}:
 [1, -1, -1]
 [-1, 1, -1]
 [-1, -1, 1]
 [1, 1, 1]
\end{minted}
The remaining four Platonic solids are the regular $3$-cube, the $3$-dimensional regular cross polytope (also known as the octahedron), the dodecahedron, and the icosahedron.
\begin{minted}{jlcon}
julia> [ f_vector(P) for P in [T, cube(3), cross_polytope(3),
                               dodecahedron(), icosahedron()] ]
5-element Vector{Vector{ZZRingElem}}:
 [4, 6, 4]
 [8, 12, 6]
 [6, 12, 8]
 [20, 30, 12]
 [12, 30, 20]
\end{minted}
The dodecahedron and icosahedron do not admit regular realizations with rational coordinates.
\OSCAR supports arbitrary real number fields with fixed embeddings.
The coordinates are exact algebraic numbers; the floating-point values inside the brackets indicate approximations.
\begin{minted}{jlcon}
julia> D = dodecahedron()
Polytope in ambient dimension 3 with EmbeddedAbsSimpleNumFieldElem type coefficients

julia> vertices(D)[1]
3-element PointVector{EmbeddedAbsSimpleNumFieldElem}:
 1//2 (0.50)
 1//4*sqrt(5) + 3//4 (1.31)
 0 (0.00)
\end{minted}
Natural generalizations of the Platonic solids vary, as does the terminology in the literature; here we follow \cite[Chapter 19]{Gruenbaum:2003}.
For instance, a polytope is \emph{semiregular} if each facet is regular and the group of symmetries acts transitively on the vertices.
It turns out that, in dimension three, the semiregular polytopes consist of the Platonic solids, the prisms and antiprisms, plus 13 exceptional types known as the \emph{Archimedean solids}.

This can be generalized further, by omitting the vertex-transitivity condition, and we arrive at \emph{regular-faced} (or \emph{uniform}) polytopes.
In dimension three this gives us the 92 types of \emph{Johnson solids}, in addition to the semiregular ones; see \cite{J66,Zalgaller:1969}.
Here the polytopes of each \emph{type} are unique, again up to scaling and rigid motions.

\begin{figure}[t]
  \centering
  \includegraphics[width=.3\linewidth]{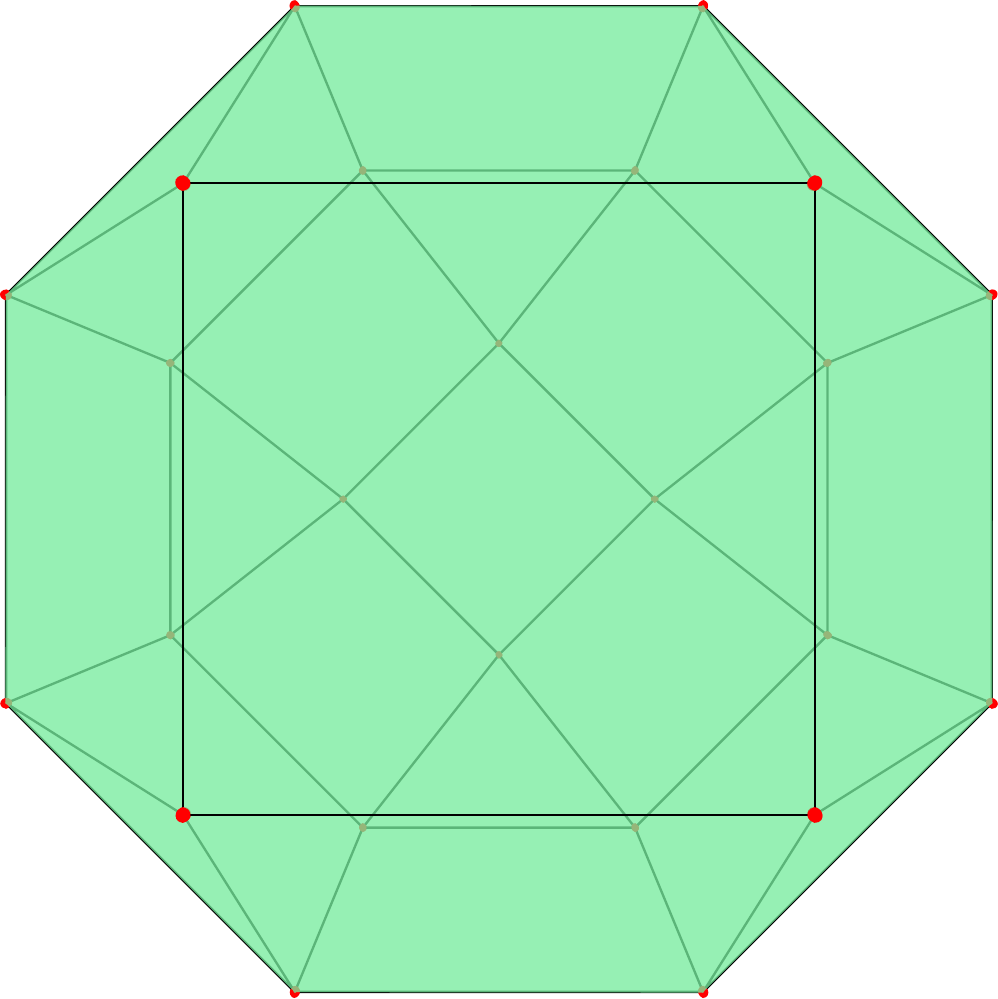}
  \caption{Elongated square gyrobicupola.  Picture generated via \mintinline{jl}{visualize(johnson_solid(37))}}
  \label{fig:J37}
\end{figure}

The \emph{elongated square gyrobicupola} is the Johnson solid $J_{37}$; it is shown in Figure~\ref{fig:J37}.
We follow the naming conventions from \cite[Table~3]{Zalgaller:1969}.
That particular polytope plays a special role, as it is sometimes counted as a 14th Archimedean solid, because any two of its vertex figures are combinatorially equivalent.
Yet it does not have a vertex transitive group of symmetries, so it is not an Archimedean solid according to  our definition.
The Johnson solid $J_{37}$ has 26 facets, eight of which are regular triangles, and the remaining 16 facets are squares.

Like the dodecahedron, the polytope $J_{37}$ does not admit a rational realization which is regular(-faced).
However, each $3$-dimensional polytope can be perturbed slightly into a rational polytope, without changing the combinatorial type.
In contrast, in higher dimensions there are polytopes which do not admit any rational realization.
A prominent example has been found by Perles (see \cite[p.94f]{Gruenbaum:2003}); in \OSCAR this is available as \mintinline{jl}{perles_nonrational_8_polytope()}.

We now turn to an example which shows how to browse and classify a set of polytopes.
The following code iterates through the 92 Johnson solids, determines the vertex orbits, and collects in a multiset how many of them have a given number of orbits. \label{GrpToPolyhedral}
\begin{minted}{jl}
orbit_counts = MSet{Int}();

for k in 1:92
  J = johnson_solid(k)
  Aut = automorphism_group(J)[:on_vertices]
  push!(orbit_counts, length(orbits(Aut)))
end
\end{minted}
\noindent
We can sort the multiset by referring to the underlying dictionary.
\begin{minted}{jlcon}
julia> sort(collect(orbit_counts.dict))
14-element Vector{Pair{Int64, Int64}}:
  2 => 26
  3 => 24
  4 => 17
  5 => 8
  7 => 5
  8 => 2
  9 => 2
 12 => 1
 14 => 1
 15 => 1
 17 => 1
 20 => 1
 27 => 1
 29 => 2
\end{minted}
The symmetry groups of the Johnson solids were discussed in \cite[\S5]{J66}.

\subsubsection{Random polytopes}\label{poly:random}
Other interesting, yet totally different, classes of polytopes come from random constructions.
To produce, e.g., the convex hull $R_{3,100}$ of 100 random points in $\RR^3$ whose coordinates are normally distributed with mean $0$ and standard deviation $1$ is as easy as calling \mintinline{jl}{convex_hull([randn(Float64, 3) for i in 1:100])} in \OSCAR.
Note that this would take the (random) coordinates of type \mintinline{jl}{Float64} and convert them to exact rational numbers.
Consequently, all subsequent computations are exact.

In practice, another class of random polytopes is quite valuable; see Section~\ref{poly:random-face-numbers} below for a case study.
These polytopes arise from taking the convex hull of points uniformly sampled from the unit sphere:
\begin{minted}{jl}
  rand_spherical_polytope(3, 100)
\end{minted}
Mathematically the situation is quite subtle.
A random point on the unit sphere in $\RR^d$ can be obtained by sampling from the $d$-variate standard normal distribution and projecting to the sphere; see \cite[\S3.4.1.E.6]{Knuth:TAOCP2} for the details.
Converting to exact rational coordinates almost surely gives a point of norm different from one.
Hence such a random point does not truly lie on $\Sph^{d-1}$.
However, producing rational random points of norm one is not difficult in itself.
We can achieve this by randomly sampling rational points in Euclidean space and employing stereographic projection, which is rational.
This is also supported in \OSCAR.
Calling
\begin{minted}{jl}
  rand_spherical_polytope(3, 100; distribution=:exact)
\end{minted}
yields a polytope in $\RR^3$ with exact rational coordinates of norm one.
The resulting distribution is not precisely uniform on the sphere, but it can be approximated arbitrarily close.
Moreover, this is more costly to compute.
So there is a tradeoff here.

\subsubsection{More examples}
As the random constructions above show, there is no shortage of polytopes.
On the contrary, there are many more polytopes than meaningful names for them.
Nonetheless, we want to show some constructions which are often useful.
Polytopes can be obtained as convex hulls of point orbits of linear groups.
The \emph{regular permutahedron} is a standard example in $\RR^4$.
\begin{minted}{jlcon}
julia> P3 = orbit_polytope([1,2,3,4], symmetric_group(4))
Polyhedron in ambient dimension 4

julia> dim(P3)
3
\end{minted}
\noindent
The input points lie in the affine hyperplane $\sum x_i=10$, so the permutahedron has codimension one.
It follows from the construction that the symmetric group acts naturally on this polytope.
Yet, the full combinatorial automorphism group, i.e., the automorphism group of the face lattice, is even larger, by a factor of two.
\begin{minted}{jlcon}
julia> describe(combinatorial_symmetries(P3))
"C2 x S4"
\end{minted}
\noindent
The additional symmetry comes from the reflection at the linear hyperplane $x_1-x_2-x_3+x_4=0$.

Polyhedral geometry in \OSCAR is not restricted to bounded polyhedra.
Given a collection $v_1,\ldots,v_k \in \RR^d$, we may define their \emph{positive hull}:
\[\pos(v_1,\ldots,v_k) = \{\lambda_1v_1+\lambda_2v_2+\cdots+\lambda_kv_k \mid \lambda_1,\ldots,\lambda_k \in \RR_{\geq 0}\} \subseteq \RR^d.\]
Subsets of $\RR^d$ arising in this way are called \emph{(polyhedral) cones}. Indeed, these sets are both cones (each is closed under positive linear combinations) and polyhedra (each is an intersection of finitely many half-spaces). As polyhedra, they have faces of various dimensions. Their one-dimensional faces are called \emph{rays}, and a cone is \emph{pointed} if it has the origin as a face.
A cone is pointed if and only if it does not contain any line through the origin.

For instance, this is the union of the $(+,+,+)$ and $(-,+,+)$ octants in $\RR^3$.
\begin{minted}{jlcon}
julia> Q = convex_hull([0 0 0], [0 1 0; 0 0 1], [1 0 0])
Polyhedron in ambient dimension 3
\end{minted}
\noindent
The constructor has three arguments, corresponding to the decomposition $Q=P+C+L$, where $P$ is a polytope, $C$ is a cone, and $L$ is a linear subspace.
Recall that any polyhedron can be written in this way, and we may even assume that the cone $C$ is pointed; e.g., see \cite[Theorem 1.2]{Ziegler:1995}.
The polyhedron $Q$ in our example does not have any vertices, but it has a unique minimal face, which is the vertical axis $x_1=x_2=0$.
\begin{minted}{jlcon}
julia> vertices(Q)
0-element SubObjectIterator{PointVector{QQFieldElem}}

julia> minimal_faces(Q)
(base_points = PointVector{QQFieldElem}[[0, 0, 0]], 
 lineality_basis = RayVector{QQFieldElem}[[1, 0, 0]])
\end{minted}
\noindent
Arbitrary affine subspaces can occur as the minimal faces of polyhedra.
Here, an affine subspace is written as the sum of a point and a linear subspace.

\begin{figure}\centering
  \includegraphics[scale=0.15]{\fd/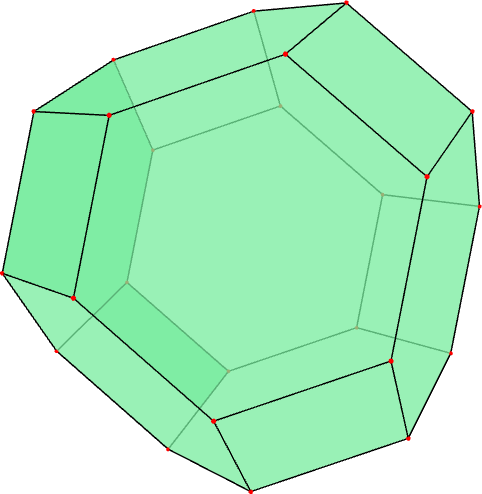} \hspace{0.3in}
  \includegraphics[scale=0.22]{\fd/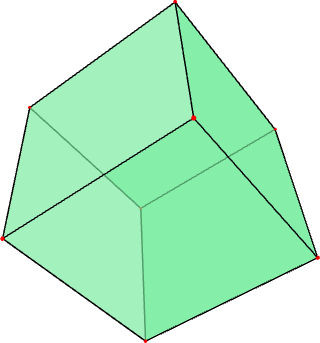}\hspace{0.3in}
  \includegraphics[scale=0.2]{\fd/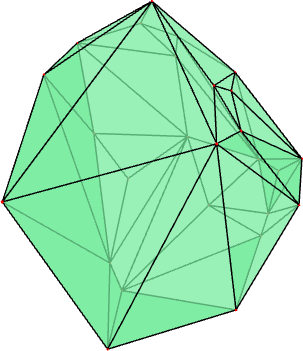}\hspace{0.3in}
  \includegraphics[scale=0.2]{\fd/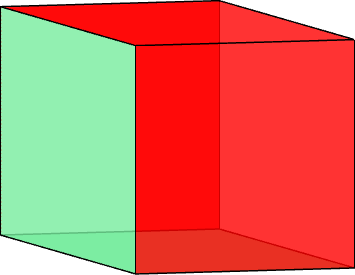}
  \caption{Visual representations of $\cP_3$, $\cC_3$, $R_{3,100}$, and $Q$. Note that $\cP_3$ is visualized within the hyperplane $x_1+x_2+x_3+x_4=10$. The polyhedron $Q$ is viewed within a box, and the unbounded portion of $Q$ is coloured red. These figures were obtained via the command \mintinline{jl}{visualize()}.} 
\end{figure}

\section{Polyhedral Fans}
Here we move on to more involved objects from polyhedral geometry.
A collection $\cF=\{C_1,\ldots,C_n\}$ of cones in $\RR^d$ is called a \emph{fan} if 
\begin{itemize}
\item Each face of each cone in $\cF$ is, itself, an element of $\cF$. 
\item For any two cones $C_1,C_2$ in $\cF$, their intersection $C_1 \cap C_2$ is also an element of $\cF$.
\end{itemize}
A fan in $\RR^d$ is \emph{complete} if the union of its cones is $\RR^d$. The fans considered in this section are usually complete, however, incomplete fans appear frequently in other important areas of mathematics; for instance, in \emph{tropical geometry} \cite{MS:2015,ETC} and \emph{toric geometry} and \cite{CLS:2011}. 

One easy way to produce a polyhedral fan in \OSCAR is by passing a collection of rays, and an incidence of rays in cones: the $(i,j)$-th entry of the incidence matrix is $1$ if cone $i$ contains ray $j$ and zero otherwise.  

\begin{minted}{jlcon}
julia> Rays = [1 0; 1 1; 0 1; -1 0; 0 -1];

julia> Cones = IncidenceMatrix([[1,2], [2,3], [3,4], [4,5], [5,1]]);

julia> PF = polyhedral_fan(Cones, Rays)
Polyhedral fan in ambient dimension 2
\end{minted}

\begin{figure}[th]\centering
\includegraphics[scale=0.25]{\fd/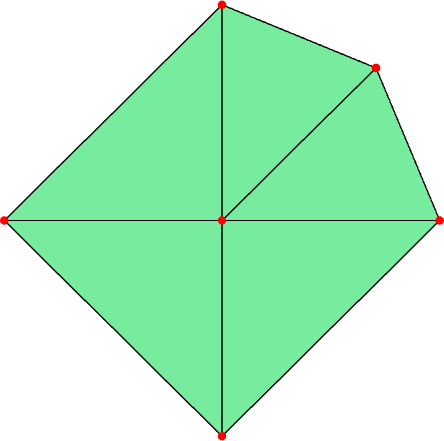}
\caption{The output of \mintinline{jl}{visualize} applied to the polyhedral fan consisting of the cones $\pos(e_1,e_1+e_2)$, $\pos(e_1+e_2,e_2)$, $\pos(e_2,-e_1)$, $\pos(-e_1,-e_2)$, $\pos(-e_2,e_1)$.}
\label{fig:PolyhedralFan}
\end{figure}

Figure \ref{fig:PolyhedralFan} shows the output of the command \mintinline{jl}{visualize(PF)}. 
Visualizing fans is not entirely trivial because fans are unbounded.
Mathematically, it would be satisfying to pass to the intersection with the unit sphere, yielding a spherical polytopal complex, which is always compact.
Yet, since spherical polytopal complexes are non-linear, this poses other challenges to the visualization.
In \OSCAR, when visualizing, we truncate each ray to unit length, and replace each cone by the convex hull of the origin and its truncated rays.
In this way, we always see an ordinary polytopal complex.

Associated to any nonempty polyhedron $P \subseteq \RR^d$ is its \emph{normal fan} $\cN(P)$ consisting of cones $N_F=\{\ell \in (\RR^d)^* \mid \ell \text{ is maximized on } F\}$ for every nonempty face $F$ of $P$.
For example, the normal fan of the polytope \[P = \conv\left( \begin{pmatrix} 0 \\ 0 \end{pmatrix} \begin{pmatrix} 2 \\ 0 \end{pmatrix} \begin{pmatrix} 2 \\ 1 \end{pmatrix} \begin{pmatrix} 1 \\ 2 \end{pmatrix} \begin{pmatrix} 0 \\ 2 \end{pmatrix} \right)\] obtained as the convex hull of five points, is exactly the fan displayed in Figure \ref{fig:PolyhedralFan}.
The normal cone $N_F$ comprises all those linear functions $\ell$ such that $\ell$ attains its maximum on $P$ at some point in the face $F$.
In this way, the normal fan of a polyhedron $P$ encodes all of the information pertaining to linear programs over $P$.
We note that the normal fan of a polytope is complete and that many different polytopes have the same normal fan.

\section{Case Studies}
The purpose of this section is to sketch some \OSCAR use-cases which occur in research on polyhedral geometry.
Some portions of the text require more substantial background in polyhedral geometry (and adjacent fields) than the previous sections.

\subsection{Face Numbers of Random Polytopes}\label{poly:random-face-numbers}
It is a major theme in geometric combinatorics to study the $f$-vectors of convex polytopes, and also of subclasses and generalizations.
A first result in this direction is Euler's equation, which says that $f_0-f_1+f_2=2$ if $f=f(P)$ is the $f$-vector of a $3$-polytope $P$.
To simplify the discussion, let us assume that $P$ is \emph{simplicial}, i.e., each proper face is a simplex.
Then a double-counting argument yields $2f_1=3f_2$, because each edge is contained in two facets, and each facet contains three edges (due to simpliciality).
Taking into account that $P$ must have at least four vertices, we see that the $f$-vectors of simplicial $3$-polytopes with $f_0=n$ vertices lie in the set
\begin{equation}\label{poly:f-vectors-3-polytopes}
  \bigSetOf{(n,\, 3n-6,\, 2n-4)\in\NN^3}{n\geq 4} \enspace ,
\end{equation}
and it is an exercise to show that for each such $f$-vector, a suitable $3$-polytope exists.
The situation in higher dimensions is much more involved.
A complete classification of the $f$-vectors of simplicial polytopes of arbitrary dimension was conjectured by McMullen \cite{McMullen:1970} and established by Stanley \cite{Stanley:1980}, Billera and Lee \cite{BilleraLee:1980}.
It is worth noting that Stanley's proof exploits the cohomological properties of the toric variety associated with a simplicial polytope in an essential way.

Now we want to test, empirically, to what extent we can explore $f$-vectors of simplicial polytopes by certain random constructions.
We start with one example: 30 vertices in six dimensions.
\begin{minted}{jlcon}
julia> P = rand_spherical_polytope(6,30)
Polytope in ambient dimension 6

julia> show(f_vector(P))
ZZRingElem[30, 336, 1468, 2874, 2568, 856]
\end{minted}
\noindent
The construction picks 30 points uniformly at random on the unit sphere $\Sph^5$ in $\RR^6$ and takes their convex hull.
We note that, almost surely, this construction produces a simplicial polytope. 
We stress that since the output above is random, reproducing it with the corresponding commands would be incredibly unlikely. 

Our goal is to see how this $f$-vector fits in the set of all possible $f$-vectors of simplicial $6$-polytopes with 30 vertices.
To this end we first convert that vector into a different coordinate system.
The \emph{$h$-vector} of the simplicial $d$-polytope $P$ is given by
\[
  h_k \ = \ \sum_{i=0}^k(-1)^{k-i}\binom{d-i}{d-k}f_{i-1} \quad \text{for } 0 \leq k \leq d \enspace,
\]
where we use the convention $f_{-1}=1$.
In particular, we have $h_0=1$, $h_1=f_0-d$ and
\[
  h_d \ = \ f_{d-1} - f_{d-2} + f_{d-3} - \dots + (-1)^{d-1}f_0 + (-1)^d  \ = \ 1 \enspace,
\]
which, up to sign, is the reduced Euler characteristic of the boundary sphere $\partial P \approx \Sph^{d-1}$.
Since the $h$-vector is obtained from the $f$-vector by an invertible linear transformation, it carries exactly the same information.
Yet the $h$-vector satisfies $h_k=h_{d-k}$, which exhibits a redundancy known as the \emph{Dehn--Sommerville equations}. The analogous symmetry in terms of the $f$-vector is more cumbersome; see \cite[\S8.3]{Ziegler:1995}.
\begin{minted}{jlcon}
julia> show(h_vector(P))
ZZRingElem[1, 24, 201, 404, 201, 24, 1]
\end{minted}
\noindent
The redundancy in the $h$-vector is eliminated by applying a final transformation 
\[
  g_0 \ = \ h_0 = 1 \quad \text{and} \quad g_k \ = \ h_k - h_{k-1} \text{ for }  1 \leq k \leq \lfloor \tfrac{d}{2}\rfloor \enspace ,
\]
to obtain the \emph{$g$-vector} of $P$.
For our example above it reads as follows.
\begin{minted}{jlcon}
julia> show(g_vector(P))
ZZRingElem[1, 23, 177, 203]
\end{minted}

\begin{figure}[th]\centering
  \includegraphics[width=.75\textwidth]{\fd/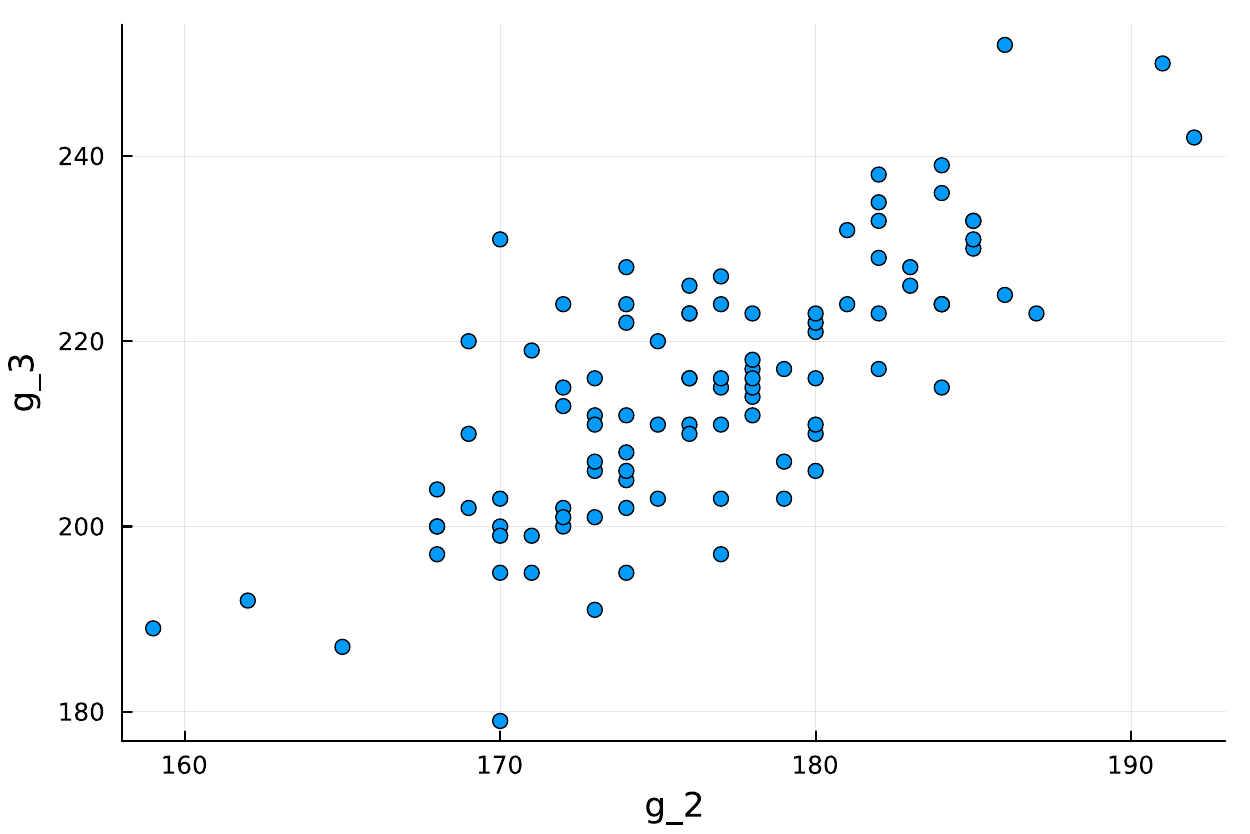}
  \caption{$g$-vectors of simplicial random $6$-polytopes with $30$ vertices.}
  \label{fig:g-vectors-scatter}
\end{figure}

Now we run a small experiment, producing $100$ random  $6$-polytopes with $30$ vertices each, and keeping track of the $g$-vectors produced.
\begin{minted}{jl}
n_vertices = 30;
n_samples = 100;
g_vectors = Array{Int32}(undef,n_samples,2);

for i=1:n_samples
    RS = rand_spherical_polytope(6,n_vertices)
    g = g_vector(RS)
    g_vectors[i,1] = g[3] # notice index shift as Julia counts from 1
    g_vectors[i,2] = g[4]
end
\end{minted}
It is instructive to plot a scatter diagram of the resulting $(g_2,g_3)$-pairs, using the standard \Julia package \mintinline{jl}{Plots}; the result is Figure~\ref{fig:g-vectors-scatter}.
\begin{minted}{jl}
using Plots
scatter(g_vectors[:,1], g_vectors[:,2],
        xlabel="g_2", ylabel="g_3", legend=false);
\end{minted}

In our sample the $g_2$-values range between 159 and 192.
This is well below the maximal possible value for $g_2$, which equals 276, by McMullen's Upper Bound Theorem.
\begin{minted}{jlcon}
julia> min_g2 = minimum(g_vectors[:,1])
159

julia> max_g2 = maximum(g_vectors[:,1])
192

julia> show(upper_bound_g_vector(6,30))
[1, 23, 276, 2300]
\end{minted}

\subsection{Gelfand--Tsetlin Polytopes}
Gelfand--Tsetlin polytopes arise in the representation theory of Lie groups. For background on the subject, see \cite{Ful91}. Given an irreducible representation $V_{\lambda}$ of $\GL(n)$ with highest weight $\lambda=(\lambda_1 \geq \ldots \geq \lambda_n)$, its restriction to a subgroup $\GL(n-1) \hookrightarrow \GL(n)$ is no longer irreducible, but decomposes as $\bigoplus_{\mu} V_{\mu}$. Here, each $V_\mu$ is an irreducible $\GL(n-1)$ representation indexed by a partition $\mu$ of $(n-1)$, and each $V_\mu$ appears with multiplicity one. The indexing set of the direct sum is given by partitions $\mu$ which \emph{interlace} $\lambda$, that is, \[\lambda_1 \geq \mu_1 \geq \lambda_2 \geq \cdots \geq \mu_{n-1} \geq \lambda_n.\] Fixing a chain of subgroups $\GL(1) \subset \GL(2) \subset \cdots \subset \GL(n)$, one may iterate this restriction process, terminating at $\textrm{dim}(V_{\lambda})$-many representations of $\GL(1)$, each of dimension one. A \emph{Gelfand--Tsetlin} diagram encodes the successive partitions indexing the representations appearing in such a chain of restrictions.

\begin{figure}[th]
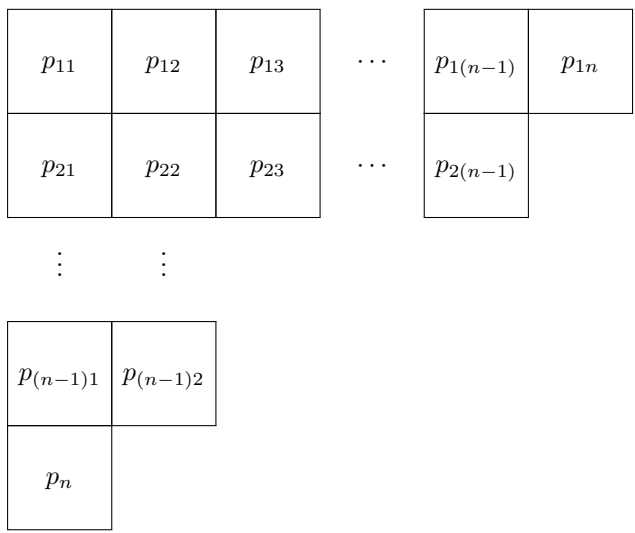

\ytableausetup{mathmode, boxsize=3.9em}
\begin{center}
\begin{ytableau}
  \none[] & p_{11} & p_{12} & p_{13} & \none[\cdots] & {p_{1(n-1)}} & p_{1n} \\
  \none[] & p_{21} & p_{22} & p_{23} &\none[\cdots]  &p_{2(n-1)}& \none \\
  \none[] & \none[\vdots] &\none[\vdots] & \none & \none &\none &\none  \\
  \none[] & p_{(n-1)1} &p_{(n-1)2} & \none & \none &\none &\none  \\
  \none[] & p_{n} & \none & \none & \none &\none &\none  \\
\end{ytableau}
\caption{A Gelfand--Tsetlin diagram.}
\label{fig:GTdiagram}
\end{center}
\end{figure}

For the combinatorics of Gelfand--Tsetlin polytopes, we follow the exposition of  \cite[\S15]{Postnikov:2009}.
A \emph{Gelfand--Tsetlin diagram} is a triangular array (see Figure \ref{fig:GTdiagram}) of ${{n+1}\choose{2}}$ real numbers $(p_{ij} \mid i,j \geq 1, i+j \leq n+1)$ 
satisfying:
\begin{itemize}
\item The first row is weakly decreasing: $p_{11} \geq p_{12} \geq \cdots \geq p_{1n}$. 
\item Each row weakly \emph{interlaces} the row above it: $p_{ij} \geq p_{(i+1)j} \geq p_{i(j+1)}$.
\end{itemize}

Given $\lambda = (\lambda_1 \geq \lambda_2 \geq \cdots \geq \lambda_n) \in \RR^n$, the \emph{Gelfand--Tsetlin polytope} $\GT(\lambda)$ is the set of all Gelfand--Tsetlin diagrams in $\RR^{{n+1}\choose{2}}$ with $\lambda$ as its first row. These are indeed polytopes, defined by inequalities of the form $p_{ij} \geq p_{(i+1)j} \geq p_{i(j+1)}$.  
When $\lambda$ is an integer partition, the lattice points of $\GT(\lambda)$ are in bijection with chains of weights of irreducible representations appearing in successive restrictions as described above. Hence, the number of lattice points in $\GT(\lambda)$ is the dimension of $V_{\lambda}$. In fact, these lattice points may be taken as a basis of $V_{\lambda}$, called a Gelfand--Tsetlin basis.

\subsubsection{An Explicit Gelfand--Tsetlin Polytope and the Weyl Dimension Formula}
\begin{figure}[hb]
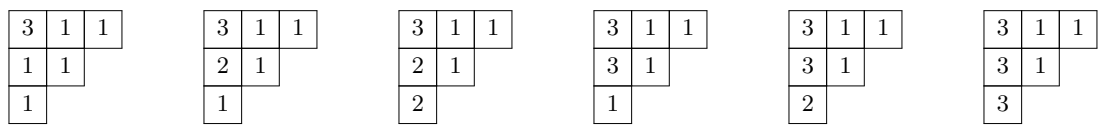

\ytableausetup{mathmode, boxsize=1.5em}
\small
\begin{ytableau}
3 & 1 & 1 \\
1 & 1 & \none \\ 
1 & \none & \none 
\end{ytableau}
\hfill
\begin{ytableau}
3 & 1 & 1 \\
2 & 1 & \none \\ 
1 & \none & \none 
\end{ytableau}
\hfill
\begin{ytableau}
3 & 1 & 1 \\
2 & 1 & \none \\ 
2 & \none & \none 
\end{ytableau}
\hfill
\begin{ytableau}
3 & 1 & 1 \\
3 & 1 & \none \\ 
1 & \none & \none 
\end{ytableau}
\hfill
\begin{ytableau}
3 & 1 & 1 \\
3 & 1 & \none \\ 
2 & \none & \none 
\end{ytableau}
\hfill
\begin{ytableau}
3 & 1 & 1 \\
3 & 1 & \none \\ 
3 & \none & \none 
\end{ytableau}
\caption{The six Gelfand--Tsetlin diagrams with integer content and $(3,1,1)$ as their top row.}
\label{fig:GTexamples}
\end{figure}

 The following code computes $\GT(3,1,1)$ and its six lattice points. Each lattice point corresponds to a diagram in Figure \ref{fig:GTexamples}. Additionally, we compute the \emph{Ehrhart polynomial} of the polytope $\GT(3,1,1)$: the univariate polynomial $\ehr(\GT(3,1,1))$ whose evaluation at $x=k$ counts the number of lattice points in the $k$-th dilate of $\GT(3,1,1)$.
\begin{minted}{jlcon}
julia> lambda = Partition([3,1,1])
[3, 1, 1]

julia> GT = gelfand_tsetlin_polytope(lambda)
Polyhedron in ambient dimension 6

julia> lattice_points(GT)
6-element SubObjectIterator{PointVector{ZZRingElem}}:
 [3, 1, 1, 1, 1, 1]
 [3, 1, 1, 2, 1, 1]
 [3, 1, 1, 2, 1, 2]
 [3, 1, 1, 3, 1, 1]
 [3, 1, 1, 3, 1, 2]
 [3, 1, 1, 3, 1, 3]

julia> ehrhart_polynomial(GT)
2*x^2 + 3*x + 1

julia> volume(project_full(GT))
2
\end{minted}

The example above shows that $\dim(V_{(3,1,1)}) = |(\GT(3,1,1) \cap \ZZ^{6})| = 6$ and more generally, that $\dim(V_{(3k,k,k)}) = \ehr(\GT(3,1,1))(k) = 2k^2+3k+1$ as shown by the Ehrhart polynomial. Alternatively, this Ehrhart computation may be verified directly from the Weyl-dimension formula 
\[\ehr(\GT(\lambda))(k) = \dim(V_{k\lambda}) = \prod_{1 \leq i,j \leq n} \frac{k\lambda_i-k\lambda_j +i-j}{j-i}.\]
From this point-of-view, the volume of $\GT(\lambda)$ is the top coefficient of $\ehr(\GT(\lambda))$: 
\begin{equation}
\label{eq:WeylDimensionVolume}
\vol(\GT(\lambda)) = \prod_{1 \leq i,j \leq n} \frac{\lambda_i-\lambda_j}{j-i}
\end{equation}
Calculating the volume of $\GT(\lambda)$ directly involves computing a triangulation of $\GT(\lambda)$, however, evaluating \eqref{eq:WeylDimensionVolume} is essentially instantaneous.

\subsubsection{Generalized Gelfand--Tsetlin Polytopes and Demazure Characters} 

In \cite{PostnikovStanley:2009}, Postnikov and Stanley introduce \emph{generalized Gelfand--Tsetlin polytopes}, which are constructed from a partition $\lambda$ as before, and a permutation $\sigma$ in the symmetric group $\Sym{n}$. 
Given a permutation $\sigma \in \Sym{n}$, its \emph{code} is the vector whose $i$-th coordinate enumerates inversions of $\sigma$ greater than $i$:
\[
\code(\sigma) = (c_1,\ldots,c_n) \quad \text{ where } \quad c_i = \bigl|\{j \mid i<j\leq n \,\, \sigma_i>\sigma_j\}\bigr|
\]
Let $w_0=(1\,n\,n-1\,\ldots,1) \in \Sym{n}$ be the long word. 
The generalized Gelfand--Tsetlin polytope $\GT(\lambda,\sigma)$ is the intersection of $\GT(\lambda)$ with the linear constraints that the first $c_i(w_0 \cdot \sigma)+1$ boxes in the $i$-th column of a Gelfand--Tsetlin diagram are equal.

For example, the permutation $w_0^2 =(1\,2\,3) \in \Sym{3}$ has $\code(w_0^2) = (1,1,0)$ as its code. Taking $\lambda=(3,1,1)$ as before, the polytope $\GT(\lambda,w_0)$ is the intersection of $\GT(\lambda)$ with the equalities $p_{11}=p_{21}$ and $p_{12}=p_{22}$. One can see explicitly that this polytope contains three lattice points, namely, the last three Gelfand--Tsetlin diagrams in Figure \ref{fig:GTexamples}: each of these diagrams has the property that the first two numbers in the first two columns are equal. 

\begin{minted}{jlcon}
julia> w0 = @perm (1,3,2)
(1,3,2)

julia> genGT = gelfand_tsetlin_polytope(lambda,w0)
Polyhedron in ambient dimension 6

julia> lattice_points(genGT)
3-element SubObjectIterator{PointVector{ZZRingElem}}:
 [3, 1, 1, 3, 1, 1]
 [3, 1, 1, 3, 1, 2]
 [3, 1, 1, 3, 1, 3]
\end{minted}

Associated to $\lambda$ and $\sigma$ is the so called \emph{Demazure module} $V_{\lambda,\sigma}$ of the complex generalized flag manifold $G/B$. Corollary 15.2 of \cite{PostnikovStanley:2009} states that when $\sigma$ is $312$-avoiding, the character $\ch_{\lambda,\sigma}(z)$ of $V_{\lambda,\sigma}$ is given by a sum over the lattice points in $\GT(\lambda,\sigma)$:
\begin{equation}
\label{eq:character}\ch_{\lambda,\sigma}(z_1,\ldots,z_n) = \sum_{P \in \GT(\lambda,\sigma) \cap \mathbb{Z}^{{n+1}\choose{2}}} z^{\weight(P)}
\end{equation}
where $\weight(P)$ is the vector which encodes the increase in the sum of the rows of $P$ (from the bottom to the top). For example, the first Gelfand--Tsetlin diagram in Figure \ref{fig:GTexamples} has weight $(1-0,2-1,5-2)=(1,1,3)$. The corresponding monomial, if this diagram were to index a summand of \eqref{eq:character}, is $z_1^1z_2^1z_3^1$.
The following code computes $\ch_{\lambda,\sigma}(z_1,z_2,z_3)$ for  $\lambda=(3,1,1)$ and $\sigma=(1,3,2)$:
\begin{minted}{jlcon}
julia> dc = demazure_character(lambda,w0)
x1^3*x2*x3 + x1^2*x2^2*x3 + x1*x2^3*x3

julia> dc(1,1,1)
3
\end{minted}
This final evaluation shows that \[\dim(V_{(3,1,1),(1,3,2)}) = |\GT((3,1,1),(1,3,2)) \cap \mathbb{Z}^6| = \ch_{(3,1,1),(1,3,2)}(1,1,1).\] Corollary 14.6 of \cite{PostnikovStanley:2009} states that this is also equal to the determinant of the matrix of binomial coefficients
\[\dim V_{\lambda,\sigma}=\det \left( {\lambda_i+n-c_i(w_0\cdot \sigma)-i}\choose{n-c_i(w_0\cdot \sigma)-j} \right)_{i,j=1}^n  \] 
which evaluates, in the case of $\lambda=(3,1,1)$ and $\sigma=w_0=(1,3,2)$ to 
\[\dim (V_{(3,1,1),(1,3,2)}) = \det \begin{pmatrix} 4 & 1 & 0 \\ 1 & 1 & 1 \\ 0 & 0 & 1\end{pmatrix} = 3. 
\]

\subsection{Hyperdeterminants and Secondary Polytopes}
Amongst all polynomials of the form 
\[F = \sum_{v\in \{0,1\}^n} c_vx_1^{v_1}x_2^{v_2} \cdots x_n^{v_n}, \quad c_v \in \CC,\]
those which have singular zero-sets form a hypersurface in the $2^n$-dimensional vector space of (complex) coefficients. This hypersurface $\cD_{2,2,\ldots,2}$ is called the ${(2 \times 2 \times \cdots \times 2)}$-hyperdeterminant. The polynomial $F \in \CC[c][x]$ is supported on the monomials whose exponent vectors $\cA$ form the vertices of the $n$-cube $\cC_n = [0,1]^3$. In other words, for nonzero choices of the coefficients $c$, the resulting polynomial has $\cC_n$ as its \emph{Newton polytope}. From this point-of-view, $\cD_{2,2,\ldots,2}$ is realized as the \emph{$\cA$-discriminant} of the $n$-cube.

A seemingly unrelated construction associated to $\cC_n$ is that of its \emph{secondary polytope}. The secondary polytope of a point configuration $\cA$ is a combinatorial object derived from information about the triangulations of $\conv(\cA)$. If $\cT$ is a triangulation of the convex hull of $\cA$ (involving only the points in $\cA$) then its \emph{GKZ-vector} is
\[\phi_{\cA}(\cT) = \sum_{a \in \cA} \sum_{T \in \cT \mid a \in T} \vol(T)e_a.\]
The \emph{secondary polytope} of $\cA$ is the convex hull of all such GKZ-vectors.

The goal of this section is to illustrate the relationship between the secondary polytope of the vertices of $\cC_n$ and the $(2 \times 2 \times \cdots \times 2)$-hyperdeterminant: as shown in \cite[Section 11]{GKZ:1994}, the secondary polytope of $\cC_n$ is the Minkowski sum of the Newton polytopes of the hyperdeterminants of the faces of $\cC_n$. Our process closely follows the exposition in \cite{HSYY:2006}.

\subsubsection{The $2 \times 2 \times 2$ Hyperdeterminant}
The polynomial $\cD_{2,2,2}$ can be easily calculated using the \mintinline{jl}{eliminate} command: 
\begin{minted}{jlcon}
julia> C = cube(3, 0, 1);

julia> R, x, c = polynomial_ring(QQ, :x=>1:3, :c=>(0:1,0:1,0:1));

julia> f = sum(prod(x[i]^Int(p[i]) for i=1:3)
               * c[(Vector{Int}(p)+[1,1,1])...] for p=lattice_points(C))
x[1]*x[2]*x[3]*c[1, 1, 1] + x[1]*x[2]*c[1, 1, 0] + 
x[1]*x[3]*c[1, 0, 1] + x[1]*c[1, 0, 0] + 
x[2]*x[3]*c[0, 1, 1] + x[2]*c[0, 1, 0] + 
x[3]*c[0, 0, 1] + c[0, 0, 0]

julia> I = ideal(R, vcat([derivative(f,t) for t = x], [f]));

julia> D222 = eliminate(I,x)[1]
c[0, 0, 0]^2*c[1, 1, 1]^2 - 2*c[0, 0, 0]*c[1, 0, 0]*c[0, 1, 1]*c[1, 1, 1] - [...] + c[1, 1, 0]^2*c[0, 0, 1]^2
\end{minted}
      
Multiplying $D_{222}$ by the determinants corresponding to all proper faces of $\cC_3$ produces a polynomial $E_{222}$ in $8$ variables called the \emph{principal $2 \times 2 \times 2$ determinant}.
Observe that we need to shift the indices by the vector $(1,1,1)$ because indices in \Julia are one-based.
\begin{minted}{jlcon}
julia> CubeFacets = [lattice_points(F) for F in faces(C,2)];

julia> VariableIndices = [[c+[1,1,1] for c=Vector{Vector{Int}}(F)] 
                           for F=CubeFacets];
 
julia> FacialDeterminants = [c[v[1]...]*c[v[4]...]-c[v[2]...]*c[v[3]...] 
                             for v=VariableIndices]
6-element Vector{QQMPolyRingElem}:
 c[0, 0, 0]*c[0, 1, 1] - c[0, 1, 0]*c[0, 0, 1]
 c[1, 0, 0]*c[1, 1, 1] - c[1, 1, 0]*c[1, 0, 1]
 c[0, 0, 0]*c[1, 0, 1] - c[1, 0, 0]*c[0, 0, 1]
 c[0, 1, 0]*c[1, 1, 1] - c[1, 1, 0]*c[0, 1, 1]
 c[0, 0, 0]*c[1, 1, 0] - c[1, 0, 0]*c[0, 1, 0]
 c[0, 0, 1]*c[1, 1, 1] - c[1, 0, 1]*c[0, 1, 1]
 
julia> E222 = D222*prod(FacialDeterminants);
\end{minted}
The Newton polytope of $E_{222}$ is a $4$-polytope $P_{222}$ in $\RR^8$ with $f$-vector $(74, 152, 100, 22)$.
\begin{minted}{jlcon}
julia> NP = newton_polytope(E222);

julia> dim(NP)
4

julia> show(f_vector(NP))
ZZRingElem[74, 152, 100, 22]
\end{minted}

\begin{figure}\centering
  \includegraphics[scale=0.4]{\fd/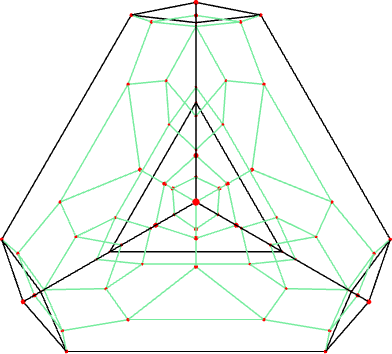} \quad \quad
  \includegraphics[scale=0.4]{\fd/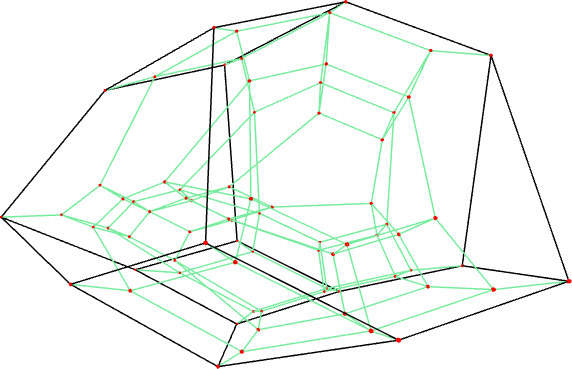}
  \caption{Two views of a Schlegel diagram for the $4$-polytope $P_{222}$ in $\RR^8$.}
  \label{fig:Schlegel}
\end{figure}

Often one chooses to illustrate a $4$-polytope like $P_{222}$ in $3$-dimensional space via a \emph{Schlegel diagram}: an image in $\RR^3$ which shows the combinatorial configuration of the $22$ facets of $P_{222}$ via a polyhedral subdivision of a $3$-polytope into $21$ polyhedra. The $22$-nd facet of $P_{222}$ is understood to correspond to the complement of this polytope. We provide two snapshots of this Schlegel diagram in Figure \ref{fig:Schlegel}.

\subsubsection{The Secondary Polytope and Triangulations of $\cC_3$}
The secondary polytope of a given polytope $P$ is interpreted in \OSCAR and \mintinline{jl}{Polymake} to be the secondary polytope associated to the vertices of $P$. Constructing such an object in \OSCAR is straightforward.
\begin{minted}{jlcon}
julia> C = cube(3,0,1)
Polytope in ambient dimension 3

julia> SP = secondary_polytope(C)
Polytope in ambient dimension 8
\end{minted}
The main connection with the previous section on hyperdeterminants is that the secondary polytope of $\cC_3$ \emph{is} the Newton polytope of $E_{222}$, up to some minor coordinate manipulations.
\begin{minted}{jlcon}
julia> V = point_matrix(vertices(NP))[:, 4:11].+1;

julia> SP == convex_hull(V)
true
\end{minted}
\noindent
The elimination which gave us the polynomial $D_{222}$ remains in the original ring, with eleven variables.
So the code above projects out the first three coordinates.
Adding the constant $1$ to each coefficient of each vertex is related to the index shift.

Alternatively, one may take the definitional path towards computing a secondary polytope by computing all triangulations of the vertices of $P$, finding their GKZ-vectors, and taking the convex hull of those. 
\begin{minted}{jlcon}
julia> AllTriangulations = all_triangulations(C);

julia> Tri_as_SOP = [subdivision_of_points(C,T)
                     for T in AllTriangulations];

julia> GKZ_Vectors = [gkz_vector(T) for T in Tri_as_SOP];

julia> SP_from_GKZ = convex_hull(GKZ_Vectors)
Polyhedron in ambient dimension 8

julia> SP_from_GKZ == SP
\end{minted}

\begin{figure}[!htpb]
  \centering
  \includegraphics[scale=0.3]{\fd/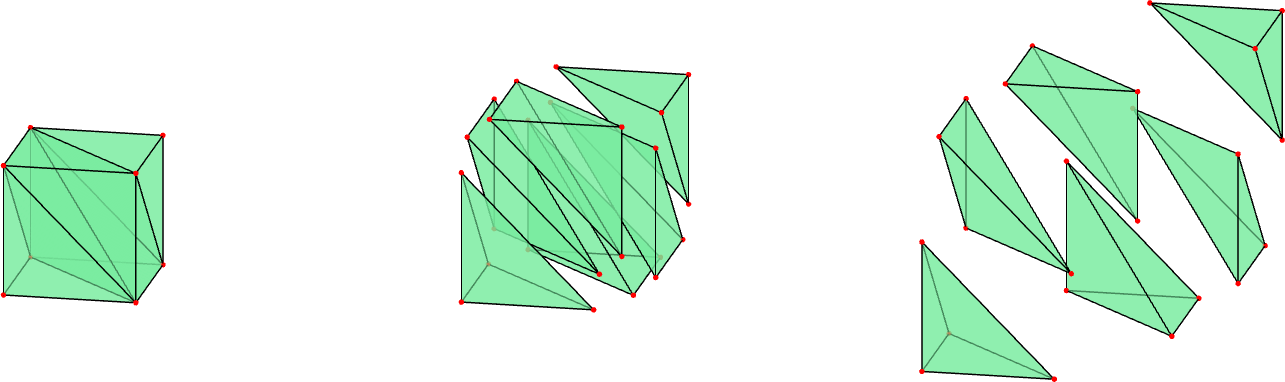}
  \caption{Three views of a triangulation of the $3$-cube with varying parameter choices for \mintinline{jl}{explode}.}
  \label{fig:Explosion}
\end{figure}

The raw output for a triangulation in \OSCAR is an array of arrays: each array at the first level corresponds to a simplex and hence contains the indices of the vertices of that simplex. For many computations which require great efficiency, such a representation is preferred. However, to obtain contextual properties of a triangulation $T$, like the GKZ-vector, one must provide the relevant polytope $P$ to \OSCAR . Warning: the ordering of the vertices of $P$ matters when constructing a subdivision of points, since the raw triangulation only refers to indices!
\begin{minted}{jlcon}

julia> T = AllTriangulations[1]
6-element Vector{Vector{Int64}}:
 [1, 2, 3, 5]
 [2, 3, 4, 5]
 [2, 4, 5, 6]
 [3, 4, 5, 7]
 [4, 5, 6, 7]
 [4, 6, 7, 8]

julia> S = subdivision_of_points(C,T)
Subdivision of points in ambient dimension 3

julia> is_regular(S)
true

julia> show(min_weights(S))
[3, 1, 1, 0, 0, 0, 0, 1]
\end{minted}
\noindent
The above evaluation checks that the triangulation $T$ is regular, that is, it is induced as the lower hull of a lifting of the vertices of $C$ by some weight function. Since it is regular, we also compute such a weight function. Finally, we showcase how to visualize a triangulation, making full-use of the visualization features \mintinline{jl}{transparency} and the famous \mintinline{jl}{explode} feature (see Figure \ref{fig:Explosion}).

\begin{minted}{jlcon}
julia> IM = IncidenceMatrix(T);

julia> V = matrix(QQ, vertices(C));

julia> PC = polyhedral_complex(IM,V);

julia> visualize(PC)
\end{minted}

The action of the symmetries of the cube induces an action on the GKZ-vectors by permutation of their coordinates. We compute the minimal representatives for the orbits these GKZ-vectors along with the number of elements in each orbit, as done in \cite{HSYY:2006}.

\begin{minted}{jlcon}
julia> G = automorphism_group(C)[:on_vertices]
Permutation group of degree 8

julia> OrbitRepresentatives=
        unique([minimum(gset(G,permuted,[v])) 
          for v in GKZ_Vectors])
6-element Vector{Vector{QQFieldElem}}:
 [1, 3, 3, 5, 5, 3, 3, 1]
 [1, 3, 3, 5, 6, 2, 2, 2]
 [1, 3, 4, 4, 6, 2, 1, 3]
 [1, 4, 4, 3, 6, 1, 1, 4]
 [1, 5, 5, 1, 5, 1, 1, 5]
 [2, 2, 2, 6, 6, 2, 2, 2]

julia> OrbitSizes = 
        [length(
          filter(x->minimum(gset(G,permuted,[x]))==u,
          GKZ_Vectors)
         ) for u in OrbitRepresentatives];

julia> show(OrbitSizes)
[12, 24, 24, 8, 2, 4]
\end{minted}

We close this section with some general remarks about computing triangulations and secondary polytopes; see also \cite{DLRS:2010}.
This is a truly demanding task and computationally infeasible for most point sets.
For instance, the 4-dimensional cube is known to have 87\,959\,448 regular triangulations, which come in 235\,277 orbits.
That is to say, the secondary polytope of $C_4$, which is 11-dimensional, has nearly 90 million vertices.
The precise number of (regular) triangulations of $C_n$ is unknown for $n\geq 5$.
State of the art software for finding all triangulations is \topcom \cite{topcom}.
For finding all \emph{regular} triangulations \mptopcom \cite{mptopcom} holds the current records.
These two packages share parts of their code bases, both run in parallel, but there are significant technical differences.
\mptopcom has been developed with large computer clusters in mind, and this makes the software not very suitable for integration into \OSCAR.
Nonetheless, finding all orbits of regular triangulations of $C_4$ with \mptopcom takes less than 10 minutes on a standard desktop computer.
On the other hand, the native implementation of secondary polytopes in \OSCAR is little more than a proof of concept.
Pfeifle gave an explicit construction of the secondary polytope of $C_3$ \cite{sec-poly-3cube}.

\section{Computing Convex Hulls}
\label{sec:ConvexHullsAndLPs}

The following two computational tasks are fundamental to polyhedral geometry.
The \emph{convex hull problem} asks to convert a finite set $V\subset\RR^d$ into finitely many affine half spaces whose intersection is the polytope $\conv(V)$.
The \emph{dual convex hull problem} is the reverse operation, which is algorithmically equivalent to its primal incarnation thanks to cone polarity.
These problems enjoy a close relationship to linear programming: computing all of the vertices of a polytope, given inequalities, is the same as solving all (generic) linear programs with these constraints.
For the geometric background of this discussion, see the textbooks \cite{Gruenbaum:2003,Ziegler:1995,JT:2013}.
In order to be computationally accessible, the coordinates are usually restricted to suitable subfields of $\RR$.
The most important example is $\QQ$, but as we saw earlier, real algebraic number fields commonly occur too.
Generally speaking, any ordered field works, e.g., ordered fields of rational functions occur at the intersection of optimization and tropical geometry \cite[\S8.5]{ETC}.

Here we collect some general wisdom about convex hull computations, and how this is reflected in \OSCAR.
Many different convex hull algorithms are known and implemented; see \cite{HDCG3:convex-hulls} for a survey.
It is important to understand, that different algorithms may exhibit very different performance, depending on the kind of input.
We illustrate this on one small example, which is not even the tip of the iceberg. Since computations are cached in the polyhedron object, organizing a fair race involves defining the \textit{same} random polytope before each command. One achieves this by using the same seed value in the random construction.

\begin{minted}{jlcon}
julia> P = rand_spherical_polytope(6, 500; seed=11)
Polytope in ambient dimension 6

julia> Polymake.prefer("ppl") do
         @time n_facets(P)
       end
205.427732 seconds (7.13 G allocations: 321.966 GiB, 0.68% gc time)
58163

julia> P = rand_spherical_polytope(6, 500; seed=11);

julia> Polymake.prefer("beneath_beyond") do
         @time n_facets(P)
       end
 61.090983 seconds (335.63 M allocations: 5.309 GiB, 0.08% gc time)
58163

julia> P = rand_spherical_polytope(6, 500; seed=11);

julia> Polymake.prefer("libnormaliz") do
         @time n_facets(P)
       end
 11.167154 seconds (19.12 M allocations: 811.301 MiB, 0.07% gc time)
58163
\end{minted}

\OSCAR's implementations and interfaces are inherited from \polymake.
Here, we use the Parma Polyhedral Library (default) \cite{BHZ08}, \polymake's own implementation \texttt{beneath\_beyond}, and libnormaliz \cite{normaliz}. Observe that, in addition to the running time, the memory consumption varies considerably. 
For a more comprehensive analysis, with more drastic examples and different winners see \cite{AvisBremnerSeidel:1997, polymake:2017}.
The timings are taken on a laptop with an Intel i7-1165G7 processor (4700MHz maximum clock) running on Manjaro Linux, with kernel 6.4.6-1, \OSCAR v0.12.2.

\section*{Acknowledgements}
Michael Joswig is supported by MaRDI (Mathematical Research
Data Initiative), funded by the Deutsche Forschungsgemeinschaft (DFG), project
number 460135501, NFDI 29/1 ``MaRDI -- Mathematische
Forschungsdateninitiative''\footnote{\url{https://www.mardi4nfdi.de/}}, and
by the SFB-TRR 195--286237555 ``Symbolic Tools in
Mathematics and their
Application''\footnote{\url{https://www.computeralgebra.de/sfb/}}.
Taylor Brysiewicz is supported by an NSERC Discovery grant (RGPIN-2023-03551).

\newpage

\bibliographystyle{plain}
\bibliography{references}

 \end{document}

%% file: ARXIVpackages.tex
\usepackage{type1cm}         
\usepackage{graphicx}        
\usepackage{multicol}        
\usepackage[bottom]{footmisc}

\usepackage{amsfonts, amsmath, amssymb}
\usepackage{xspace}          
\usepackage{ytableau}        
\usepackage{pgf}

\usepackage[draft=true]{minted} 
\setminted{breaklines=true}
\definecolor{lbcolor}{rgb}{0.9,0.9,0.9}
\setminted{bgcolor=lbcolor}
\setminted{fontsize=\footnotesize}

\usepackage{mathtools}
\usepackage{booktabs}        
\usepackage{tikz}            
\usepackage{csquotes}        

\usepackage[colorinlistoftodos, bordercolor=orange, backgroundcolor=orange!20, linecolor=orange, textsize=scriptsize]{todonotes} 

\usepackage[noend]{algorithmic}

\usepackage{calc,ifthen,alltt,tabularx,capt-of}
\usepackage[linesnumbered,longend,ruled,vlined]{algorithm2e}

\usepackage[bookmarks=false]{hyperref}

\usepackage{longtable}
\usepackage{makecell}
\usepackage{wrapfig}

\usepackage[capitalise, noabbrev]{cleveref} 
\usepackage{enumitem} 

%% file: ARXIVdefinitions.tex
\newcommand\OSCAR{\texttt{OSCAR}\xspace}
\newcommand\Julia{\texttt{Julia}\xspace}

\newcommand\polymake{\texttt{polymake}\xspace}

\newcommand\topcom{\texttt{TOPCOM}\xspace}
\newcommand\mptopcom{\texttt{mptopcom}\xspace}

\newcommand\FF{\mathbb{F}}
\newcommand\NN{\mathbb{N}} 
\newcommand\ZZ{\mathbb{Z}} 
\newcommand\QQ{\mathbb{Q}} 

\newcommand\RR{\mathbb{R}} 
\newcommand\CC{\mathbb{C}} 

\newcommand\Sph{\mathbb{S}} 

\newcommand\cA{{\mathcal A}}

\newcommand\cC{{\mathcal C}}
\newcommand\cD{{\mathcal D}}

\newcommand\cF{{\mathcal F}}

\newcommand\cN{{\mathcal N}}
\newcommand\cT{{\mathcal T}}
\newcommand\cP{{\mathcal P}}

\newcommand\smallSetOf[2]{\{#1 \mid #2\}}
\newcommand\bigSetOf[2]{\bigl\{#1 \bigm| #2\bigr\}}

\DeclareMathOperator{\GL}{GL} 
\newcommand{\Sym}[1]{\mathfrak S_{#1}} 
\DeclareMathOperator{\code}{code} 

\DeclareMathOperator{\weight}{weight} 
\DeclareMathOperator{\ch}{ch} 

\DeclareMathOperator{\conv}{conv} 
\DeclareMathOperator{\pos}{pos} 
\DeclareMathOperator{\ehr}{ehr} 
\DeclareMathOperator{\GT}{GT} 
\DeclareMathOperator{\vol}{vol} 

\definecolor{promptColor}{rgb}{0.0,0.0,0.589}
\definecolor{brkpromptColor}{rgb}{0.589,0.0,0.0}
\definecolor{gapinputColor}{rgb}{0.589,0.0,0.0}
\definecolor{gapoutputColor}{rgb}{0.0,0.0,0.0}
\definecolor{darkgreen}{rgb}{0.05,0.6,0.1}
\definecolor{colrem}{rgb}{0,0.7,0}

\newcommand{\tmfloatcontents}{}
\newlength{\tmfloatwidth}
\newcommand{\tmfloat}[5]{
	\renewcommand{\tmfloatcontents}{#4}
	\setlength{\tmfloatwidth}{\widthof{\tmfloatcontents}+1in}
	\ifthenelse{\equal{#2}{small}}
	{\setlength{\tmfloatwidth}{0.45\linewidth}}
	{\setlength{\tmfloatwidth}{\linewidth}}
	\begin{minipage}[#1]{\tmfloatwidth}
		\begin{center}
			\tmfloatcontents
			\captionof{#3}{#5}
		\end{center}
\end{minipage}}

